\documentclass[12pt]{amsart}
\usepackage{amsmath,amssymb, amsfonts,amscd,psfrag,graphicx, enumitem}
\usepackage{mathrsfs}

\theoremstyle{plain}
\newtheorem{theorem}{Theorem}
\newtheorem{lemma}[theorem]{Lemma}
\newtheorem{corollary}[theorem]{Corollary}
\newtheorem{proposition}[theorem]{Proposition}

\theoremstyle{definition}
\newtheorem{definition}[theorem]{Definition}

\newtheorem{problem}[theorem]{Problem}
\newtheorem{question}[theorem]{Question}

\theoremstyle{remark}
\newtheorem{remark}[theorem]{Remark}

\newcommand{\PC}{\mathscr{PC}}

\newcommand{\N}{\mathbb{N}}

\newcommand{\F}{\mathcal{F}}

\newcommand{\IE}{\mathcal{I}}
\newcommand{\HE}{\mathcal{H}}

\newcommand{\cd}{\operatorname{cochord}}
\newcommand{\reg}{\operatorname{reg}}
\newcommand{\di}{\operatorname{diam}}
\newcommand{\gi}{\operatorname{girth}}

\newcommand{\im}{\operatorname{im}}
\newcommand{\m}{\operatorname{m}}
\newcommand{\ve}{\overrightarrow}

\newcommand{\cov}{\textrm{Cov}}

\begin{document}

\bibliographystyle{plain}

\title[Vertex-decomposability, codismantlability and regularity]{Vertex-decomposable graphs, codismantlability,
Cohen-Macaulayness, and Castelnuovo-Mumford regularity}
\author{T\" urker B\i y\i ko\u glu and Yusuf Civan}

\address{Department of Mathematics, Izmir Institute of Technology,
Gulbahce, 35437, Urla, Izmir, Turkey}

\address{Department of Mathematics, Suleyman Demirel University,
Isparta, 32260, Turkey.}

\email{turkerbiyikoglu@iyte.edu.tr\\
yusufcivan@sdu.edu.tr}


\date{\today}

\thanks{Both authors are supported by T\" UBA through Young Scientist Award
Program (T\" UBA-GEB\. IP/2009/06 and 2008/08) and by T\" UB\. ITAK, grant no:111T704}


\begin{abstract}
We call a vertex $x$ of a graph $G=(V,E)$ a \emph{codominated vertex} if $N_G[y]\subseteq N_G[x]$ for some vertex 
$y\in V\backslash \{x\}$, and
a graph $G$ is called \emph{codismantlable} if either it is an edgeless graph or it contains a codominated vertex $x$ such that
$G-x$ is codismantlable. We show that $(C_4,C_5)$-free vertex-decomposable graphs are codismantlable, and
prove that if $G$ is a $(C_4,C_5,C_7)$-free well-covered graph, then vertex-decomposability, codismantlability
and Cohen-Macaulayness for $G$ are all equivalent. These results complement and unify many of the earlier results on bipartite, 
chordal and very well-covered graphs.
We also study the Castelnuovo-Mumford regularity $\reg(G)$ of such graphs, and show that $\reg(G)=\im(G)$
whenever $G$ is a $(C_4,C_5)$-free vertex-decomposable graph, 
where $\im(G)$ is the induced matching number of $G$. Furthermore, we prove that 
$H$ must be a codismantlable graph if $\im(H)=\reg(H)=\m(H)$, where $\m(H)$
is the matching number of $H$. 
We further describe an operation on digraphs that creates a vertex-decomposable and codismantlable graph from any acyclic digraph. 
By way of application, we provide an infinite family $H_n$ ($n\geq 4$)
of sequentially Cohen-Macaulay
graphs whose vertex cover numbers are half of their orders, while containing no vertex of degree-one such that they are 
vertex-decomposable, and
$\reg(H_n)=\im(H_n)$ if $n\geq 6$. This answers a recent question of Mahmoudi, et al~\cite{MMCRTY}.
\end{abstract} 

\maketitle
\section{Introduction}
The present work is devoted to the study of algebraic and combinatorial properties of a new graph class, codismantlable graphs. 
Our results complement and unify many of the earlier results on the topic~\cite{DE, FVT1, MMCRTY, VT, RW1,RW2}. 
We prove that if a graph $G$ lacks certain induced cycles, the vertex-decomposability and (sequentially) Cohen-Macaulayness
of $G$ relies upon the codismantlability of $G$. Moreover, such an approach also permits to read off the Castelnuovo-Mumford
regularity of such graphs combinatorially.

All graphs we consider are finite and simple. Given a graph $G=(V,E)$ with vertex set $V=\{u_1,\ldots,u_n\}$.
The \emph{edge ideal} $I(G)$ of $G$ is defined to be the ideal $I(G)\subseteq R=k[x_1,\ldots,x_n]$ (where $k$ is a field)
generated by all monomials $x_ix_j$ such that $u_iu_j\in E$. 
Such an assignment interrelates algebraic properties and invariants of $I(G)$ 
to combinatorial properties and invariants of $G$, and vice versa. Most recent research in this area focuses on  
finding combinatorial properties of graphs that guarantee (sequentially) Cohen-Macaulayness of $R/I(G)$, 
as well as enable us to combinatorially compute algebraic invariants of $R/I(G)$, 
such as Castelnuovo-Mumford regularity. We refer readers to \cite{MV} for an excellent recent survey on the subject. 

Topological combinatorics enters the picture via the Stanley-Reisner correspondence. Recall that the
\emph{independence complex} $\IE(G)$ of a graph $G$ is the simplicial complex on $V$ whose faces are independent sets
of $G$. In this context, $R/I(G)$ is the Stanley-Reisner ring of $\IE(G)$. This assignment allows us to use methods from 
topological combinatorics, such as shellability and vertex-decomposability
created to analyze the structure of Stanley-Reisner rings of simplicial
complexes not necessarily those arisen from graphs. It seems that such methods are well-adapted to the study of edge rings,
and almost all existing characterizations of graphs with (sequentially) Cohen-Macaulay edge rings involve these notions. 

It has been shown that the vertex-decomposability
of the independence complex $\IE(G)$ of a graph $G$ and (sequentially) Cohen-Macaulayness of $R/I(G)$ are equivalent if the graph 
$G$ is a member of the family  of bipartite graphs~\cite{VT} or chordal graphs~\cite{DE, FVT1,RW1} or very well-covered 
graphs~\cite{MMCRTY}.  A detailed look at the proofs in these papers reveals that the vertex-decomposability of these graphs 
follows from the existence of a shedding order with a common property. Indeed, 
if a vertex-decomposable graph belongs to one of these graph classes 
with at least one edge, then it has a shedding vertex $x$ satisfying that the closed neighborhood set of one of its neighbors is
contained in that of $x$. Such an observation naturally raises a simple question: What conditions on a given graph  
guarantee that its shedding vertices (if any) satisfy this property? We obtain a fairly comprehensive answer to this question by proving 
that if a graph lacks induced cycles of length four and five, any of its shedding vertices must be of this type. 
This leads us to introduce a recursively-defined graph class, that of \emph{codismantlable graphs} 
(see Definition~\ref{defn:codis}, Corollary~\ref{cor:ver-dec-codis} and Theorem~\ref{thm:wc-cns}).

In most cases, our method of proof is of a combinatorial nature. In Section $2$, we introduce codismantlable graphs, and
reveal some of their combinatorial properties.
We prove that when a well-covered graph $G$ is $(C_4,C_5,C_7)$-free, then its vertex-decomposability, 
codismantlability and Cohen-Macaulayness
are all equivalent (see Theorem~\ref{thm:wc-cns}). In a similar vein, we show that if $G$ is a well-covered graph with $\gi(G)\geq 5$,
then $G$ is vertex-decomposable if and only if it is Cohen-Macaulay.

Section $3$ deals with the Castelnuovo-Mumford regularity of codismantlable graphs. We prove that 
the regularity of a $(C_4,C_5)$-free vertex-decomposable graph is
equal to its induced matching number. Furthermore, we show that the cochordal cover number of a graph of 
girth at least five is equal to its edge-domination number.

In the final section, we associate a (simple) graph to any given digraph, the \emph{common-enemy graph},
and prove that common-enemy graphs are always vertex-decomposable and codismantlable provided that the digraph we
begin with is acyclic. 
We show that the upper bound graphs of (finite) posets are examples of common-enemy graphs so that they are vertex decomposable
and codismantlable. As an application, we introduce the {\it edge-clique-whiskering} of a given graph with respect to any of 
its edge-clique partition, and provide an infinite family of sequentially Cohen-Macaulay
graphs whose vertex cover numbers are half of their orders, while containing no vertex of degree-one. 

\subsection{Preliminaries}
By a graph $G=(V,E)$, we will mean an undirected graph without loops or
multiple edges. An edge between $u$ and $v$ is denoted by $e=uv$ or
$e=\{u,v\}$ interchangeably. A graph $G=(V,E)$ is called an \emph{edgeless graph} on $V$ whenever
$E=\emptyset$. If $U\subset V$, the graph induced on $U$ is written $G[U]$, and in particular,
we abbreviate $G[V\backslash U]$ to $G-U$, and write $G-x$ whenever $U=\{x\}$. 
Throughout, we denote by $C_n$, the $n$-\emph{cycle}, which is the graph on $\{x_1,x_2,\ldots,x_n\}$
with edges $\{x_1x_2,\ldots,x_{n-1}x_n,x_nx_1\}$.

For a given subset $U\subseteq V$, the (open) neighborhood of $U$ is
defined by $N_G(U):=\cup_{u\in U}N_G(u)$, where $N_G(u):=\{v\in V\colon uv\in E\}$,
and similarly, $N_G[U]:=N_G(U)\cup U$ is the 
closed neighborhood of $U$. We call a vertex $x\in V$ a \emph{full-vertex} whenever
$N_G[x]=V$.

We say that $G$ is $H$-free if no induced subgraph of $G$ is isomorphic to $H$.
A graph $G$ is called \emph{chordal} if it is $C_k$-free for any $k\geq 4$, and a graph $H$ is
said to be \emph{cochordal} if its complement is a chordal graph.
A subset $S\subseteq V$ is called a {\it clique} of $G$ if $G[S]$ is isomorphic to a complete
graph. 

A graph $G$ is said to be \emph{well-covered} (or equivalently \emph{unmixed})
if every maximal independent set has the same size,
which is equivalent to requiring that $\IE(G)$ is a pure simplicial complex.
Finally, a subset $X\subseteq V$ is called a 
\emph{vertex cover} of $G$ if $e\cap X\neq \emptyset$ for any $e\in E$, and a vertex cover is
\emph{minimal} if no proper subset of it is a vertex cover for $G$. We should note that a graph $G$ is
well-covered if and only if every minimal vertex cover has the same size. 

Recall that a subset $M\subseteq E$ is called a {\it matching} of $G$ if no two edges in $M$ share a common end,
and a maximum matching is a matching that contains the largest possible number of edges. The {\it matching number}
$\m(G)$ of $G$ is the cardinality of a maximum matching. Moreover, a matching $M$ of $G$ is an {\it induced matching}
if it occurs as an induced subgraph of $G$, and the cardinality of a maximum induced matching is called the
{\it induced matching number} of $G$ and denoted by $\im(G)$. 

A vertex $x\in V$ is called a \emph{leaf} if $x$ is of degree $1$, and a \emph{pendant edge}
in $G$ is an edge which is incident to a leaf.
An edge $e$ is called a \emph{triangle edge} if its two endpoints
have degree $2$ and have a common neighbor.

When $I$ is an independent set of a graph $G$, the \emph{link} of $I$ in $\IE(G)$ is the independence
complex of $G-N_G[I]$.
A graph $G$ is said to be (sequentially) Cohen-Macaulay whenever $R/I(G)$ is. We recall that following
Reisner's theorem~\cite{RS}, if $G$ is a (sequentially) Cohen-Macaulay graph, and if $I$ is an independent set in $G$,
then the graph $G-N_G[I]$ is also a (sequentially) Cohen-Macaulay graph. In other words, the links in a (sequentially) Cohen-Macaulay
graph are (sequentially) Cohen-Macaulay.

\section{Vertex-decomposable graphs and codismantlability}

In this section, we introduce our main objects of study, codismantlable graphs, investigate some of their combinatorial properties
and determine their role in the characterization of vertex-decomposable graphs without any induced cycles of length four and five.
\begin{definition}
A vertex $x$ of $G=(V,E)$ is called {\it codominated} if there exists a vertex $y\in V\backslash \{x\}$ 
such that $N_G[y]\subseteq N_G[x]$.
\end{definition}

\begin{definition}\label{defn:codis}
Given two graphs $G$ and $H$. We say that $G$ is {\it codismantlable to} $H$ if there exist graphs $G_0,G_1,\ldots,G_{k+1}$
satisfying $G\cong G_0$, $H\cong G_{k+1}$ and $G_{i+1}=G_i-x_i$ for each $0\leq i\leq k$, where $x_i$ is codominated in $G_i$.
A graph $G$ is called {\it codismantlable} if either it is an edgeless graph or it is codismantlable to an edgeless graph.
When $G$ is a codismantlable graph, the ordered list $\{x_0,x_1,\ldots,x_k\}$ of vertices is called a \emph{codismantling order}
(or shortly a \emph{cd-order}) for $G$.
\end{definition}

\begin{remark}
We note that there is no direct relation between the classes of codismantlable and dismantlable graphs, so our terminology
is just a coincidence stemming from the relevance of their defining relations. 
We recall that a vertex $u$ of a graph $G$ is said to be \emph{dominated} by a vertex $v\in V\backslash \{u\}$
if $N_G[u]\subseteq N_G[v]$, and a graph $G$ is called a \emph{dismantlable graph} (also known as a \emph{cop-win graph}) 
if either it consists of a single vertex or it has a dominated vertex
$u$ such that $G-u$ is dismantlable \cite{BLS}. However, there are dismantlable graphs that are not codismantlable, and vice versa. 
For instance, the wheel graph $W_n$ for $n\geq 4$ is an example of a dismantlable graph that is not codismantlable, and the pan graph 
$\textrm{Pan}_m$
for $m\geq 4$ is a codismantlable graph which is not dismantlable. Moreover,
they are not complementary graph classes, that is, the complement of a codismantlable graph need not be a dismantlable graph and 
vice versa. 
\end{remark}

We recall that a vertex $v$ of a graph $G$ is called a \emph{simplicial vertex} if $N_G(x)$ is a clique of $G$.
This means that any neighbor of a simplicial vertex is a codominated vertex of the graph. So,
the following is obvious.

\begin{corollary}
Any  chordal graph is codismantlable. 
\end{corollary}

\begin{definition}
A vertex $x$ of $G$ is called a {\it shedding vertex} if for every independent set $S$ in $G-N_G[x]$, there is some vertex $v\in N_G(x)$
so that $S\cup \{v\}$ is independent. A graph $G$ is called {\it vertex-decomposable} if either it is an edgeless graph or it 
has a shedding vertex $x$ such that $G-x$ and $G-N_G[x]$ are both vertex-decomposable.
\end{definition}

We remark that the definition of vertex-decomposability originates in a statement about the independence
complex rather than the graph itself. However, this way of defining it is easily seen to be equivalent to more conventional
approaches~\cite{RW1}.

In many examples of vertex-decomposable graphs~\cite{DE, MMCRTY, VT, RW1}, the existing shedding vertices are of particular kind.

\begin{lemma}\label{lem:RW1}
If $x$ is codominated in $G$, then it is a shedding vertex.
\end{lemma}
\begin{proof}
Let $y\in V\backslash \{x\}$ be a vertex such that $N_G[y]\subseteq N_G[x]$. Then, if $S$ is an independent set
in $G-N_G[x]$, so is $S\cup \{y\}$ in $G$.
\end{proof}

Under the assumption that a graph lacks some particular induced cycles, we prove that the converse of Lemma~\ref{lem:RW1} also
holds. 

\begin{theorem}\label{thm:shedding}
Let $G$ be a $(C_4,C_5)$-free graph. Then a vertex $x$ of $G$ is a shedding vertex if and only if it is codominated in $G$.
\end{theorem}
\begin{proof}
Following Lemma~\ref{lem:RW1}, we only need to verify the necessity of the claim. So, assume that  $G=(V,E)$ is a $(C_4,C_5)$-free graph
with $x$ a shedding vertex of $G$. Suppose that $x$ is not codominated in $G$. Since $x$ is a shedding vertex, we must have
$N_G(x)\neq \emptyset$. We let $K:=N_G(x)$ and $S:=N_G(K)\backslash K$.

We notice that it suffices to find an independent set $U\subseteq S$
such that every vertex in $K$ is adjacent to a vertex of $U$.
Assume that this is not the case, and let $I\subseteq S$ be an independent set which is maximal
with the property that the number of its neighbors in $K$ is
of maximum size. Then there exists a vertex $v\in K$ such that $va\notin E$ for any $a\in I$. However, such a vertex $v$ must have
a neighbor, say $b$ in $S\backslash I$, since $x$ is not codominated in $G$. It then follows that $bc\in E$ for some $c\in I$,
since otherwise the independent set $I\cup \{b\}$ would have more neighbors in $K$ then $I$. 
If $cw\in E$ for $w\in K$, then we must have $bw\in E$. This is due to fact that if $bw\notin E$, then either
the set $\{v,b,c,w\}$ induces a $C_4$ or the set $\{x,v,b,c,w\}$ induces a $C_5$ in $G$, any of which is impossible.
Therefore, we may replace the neighbors of $b$ in $I$ with $b$ and get a new independent set contained in $S$
that is adjacent with more vertices in $K$. This contradicts the maximality of $I$.

Now, let $U\subseteq S$ be an independent set such that every vertex in $K$ is adjacent to a vertex in $U$. But,
then $U\cup \{z\}$ can not be an independent set for any $z\in K=N_G(x)$. This is a contradiction that $x$ is
a shedding vertex.
\end{proof}

If $x$ is a shedding vertex for which the set of its neighbors forms an independent set in $G$, we may further weaken the conditions of 
Theorem~\ref{thm:shedding} to conclude that $x$ is a codominated vertex, whose proof is almost similar to that of 
Theorem~\ref{thm:shedding}, so it is omitted. 

\begin{proposition}\label{prop:bip-vd}
If $G$ is $C_5$-free and $x$ is a shedding vertex such that $N_G(x)$ is an independent set, then $x$ is codominated in $G$.
\end{proposition}

\begin{corollary}\label{cor:ver-dec-codis}
Any $(C_4,C_5)$-free vertex-decomposable graph is codismantlable.
\end{corollary}

\begin{figure}[ht]
\begin{center}
\psfrag{1}{$x_0$}\psfrag{2}{$x_1$}\psfrag{3}{$x_2$}\psfrag{4}{$x_3$}
\includegraphics[width=3in,height=1in]{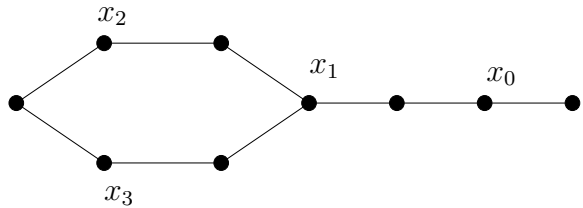}
\end{center}
\caption{A $(C_4,C_5)$-free codismantlable graph $G$ which is not vertex-decomposable.}
\label{contr-1}
\end{figure}

Clearly, the converse of Corollary~\ref{cor:ver-dec-codis} does not hold in general. The graph $G$ in Figure~\ref{contr-1}
is codismantlable with a cd-order $\{x_0,x_1,x_2,x_3\}$, while it is not vertex-decomposable, its only shedding vertex
is $x_0$, and $G-N_G[x_0]\cong C_6$ is not vertex-decomposable~\cite{FVT1}.

At this point, it is natural to seek extra conditions on graphs for which their existence may guarantee that
the converse of Corollary~\ref{cor:ver-dec-codis} holds. Having this in mind,
we first need some technical results related to well-covered graphs.

The characterization of well-covered graphs heavily depends on the existence of vertices with some extra properties~\cite{FHN}.
We recall that a vertex $x$ of a well-covered graph $G$ is called an \emph{extendable vertex} if $G-x$ is well-covered and
$\alpha(G)=\alpha(G-x)$,
which is equivalent to requiring that there exists no independent set
$I$ in $G-N_G[x]$ such that $x$ is an isolated vertex of $G-N_G[I]$  (see Lemma $2$ of~\cite{FHN}). 
It then follows that when a graph $G$ is well-covered,
the notions of shedding and extendable vertices coincide in $G$. 

\begin{lemma}\label{lem:unmixed}
If $G$ is a well-covered graph and $I$ is an independent set in $G$, then $G-N_G[I]$ is well-covered. Furthermore,
if $x$ is a codominated vertex of $G$, then $G-x$ is well-covered.
\end{lemma}
\begin{proof}
The first statement is well-known (see~\cite{FHN}), and the second follows from the fact that any such vertex is
necessarily a shedding vertex~\cite{PB}.
\end{proof}
\begin{corollary}\label{cor:unmixed}
Let $G$ and $H$ be two graphs such that $G$ is codismantlable to $H$.
If $G$ is well-covered, then so is $H$. 
\end{corollary}

\begin{corollary}\label{cor:exten-codom}
Let $G$ be a $(C_4,C_5)$-free well-covered graph. Then a vertex  $x$ of $G$ is codominated if and only if 
there exists no independent set $I$ in $G-N_G[x]$ such that $x$ is an isolated vertex of $G-N_G[I]$.
\end{corollary}
\begin{proof}
Let $N_G[y]\subseteq N_G[x]$ for some vertex $y\in V\backslash \{x\}$. If $I$ is an independent set in $G-N_G[x]$,
then $y\notin N_G[I]$. Thus, $x$ can not be isolated in $G-N_G[I]$.

The converse follows from Theorem~\ref{thm:shedding}, since such a vertex is necessarily a shedding vertex of $G$,
so it must be codominated.
\end{proof}

\begin{lemma}\label{lem:lnk-codis}
Let $G=(V,E)$ be a $(C_4,C_5,C_7)$-free well-covered graph. If $x$ is a vertex of $G$ such that $G-N_G[x]$ is codismantlable,
then $G$ is codismantlable.
\end{lemma}
\begin{proof}
It suffices by induction to find a codominated vertex in $N_G(x)$. If $N_G[x]$ is a clique, then every neighbor of $x$
is codominated. Otherwise, suppose that there is no edge between some pair of vertices $p$ and $q$ in $N_G(x)$. We consider the following
three cases.

{\it Case 1:} Neither $p$ nor $q$ has a neighbor in $H:=G-N_G[x]$. This can not happen, since if $I$ is
a maximal independent set containing $p$ and $q$, then $(I\backslash \{p,q\})\cup \{x\}$ is also a maximal,
contradicting the well-covered condition.  

{\it Case 2:} Exactly one of $p$ and $q$ has a neighbor in $H$. Assume without loss of generality that $p$
has a neighbor in $H$, while $q$ does not. We show that $p$ is codominated. Suppose to the contrary that $p$ is
not a codominated vertex. Let $y\in V(H)$ be a neighbor of $p$. It follows that $N_G[y]\nsubseteq N_G[p]$.
Note that for such a vertex $y\in V(H)$, if $yw\in E$ for some $w\in N_G(x)$, then $wp\in E$, 
since otherwise the set $\{u,c,w,x\}$ induces a $C_4$ in $G$, which is not possible. 
Therefore, there exists $v\in N_H[y]$ such that $vp\notin E$.
We define $D_p:=N_H(p)=\{s\in V(H)\colon ps\in E\}$, the set of neighbors of $p$ in $H$, and 
let $S_p:=N_H(D_p)\backslash D_p$. Since $G$ is $(C_4,C_5)$-free and $p$ is not codominated, 
we conclude by Corollary~\ref{cor:exten-codom} that 
there exists an independent set, say $U_p\subseteq S_p$ such that every vertex in $D_p$ is adjacent to a vertex in $U_p$. 
It then follows that the graph $L_p:=G-N_G[U_p]$ must be well-covered by Lemma~\ref{lem:unmixed}. 
On the other hand, this implies that $L_p-N_{L_p}[x]$ is well-covered too. However, if $I_p$ is a maximal independent of 
$L_p-N_{L_p}[x]$, then $I_p\cup \{x\}$ and $I_p\cup \{p,q\}$ are both independent sets in $L_p$, while $I_p\cup \{x\}$ being maximal.
This is a contradiction to $L_p$ being well-covered. Therefore, the vertex $p$ must be codominated,

{\it Case 3:} $p$ and $q$ both have a neighbor in $H$. We show that either $p$ or $q$ is codominated. Assume otherwise
that neither $p$ nor $q$ is codominated. If we define $D_p:=N_H(p)$ and $D_q:=N_H(q)$,
then the sets $S_p:=N_H(D_p)\backslash D_p$ and $S_q:=N_H(D_q)\backslash D_q$ are both non-empty.
Note also that if $u\in D_p$ and $v\in D_q$, then $uv\notin E$, since $G$ is $(C_4,C_5)$-free. 
Then by Corollary~\ref{cor:exten-codom}, there exist independent sets $U_p\subseteq S_p$
and $U_q\subseteq S_q$ such that $p$ and $q$ have no neighbors in $H-N_H[U_p]$ and
$H-N_H[U_q]$ respectively. Furthermore, since $G$ is $(C_4,C_5,C_7)$-free, if $t\in U_p$ and
$r\in U_q$ such that $tr\in E$, then  either $tu\in E$ for any $u\in D_q$ with $ru\in E$ or $rv\in E$ for any $v\in D_p$
with $tv\in E$. It then follows that there 
exists an independent set $U\subseteq U_p\cup U_q$ such that neither $p$ nor $q$ have any neighbor
in $H-N_H[U]$. If we set $L:=G-N_{G}[U]$, the graph $L$ is well-covered. However,
if $I$ is a maximal independent set in $L-N_L[x]$, then $I\cup \{p,q\}$ and $I\cup \{x\}$ are both independent sets of
$L$, while $I\cup \{x\}$ being maximal; hence, $L$ can not be well-covered, a contradiction. Therefore, at least one of
$p$ and $q$ must be codominated. This completes the proof.
\end{proof}
We remark that if a graph $G$ is not $C_7$-free, Lemma~\ref{lem:lnk-codis} may not need to hold in general,
as $C_7$ is itself a $(C_4,C_5)$-free well-covered graph where $C_7-N_{C_7}[x]$ is codismantlable for any $x\in V(C_7)$. 

We call a graph $G$ \emph{closed neighborhood Sperner} (or shortly a \emph{CNS-graph}) if $G$ contains no
codominated vertex. It follows that if $G$ is a CNS-graph, 
the set $\{N_G[x]\colon x\in V(G)\}$ is an anti-chain with respect to the inclusion order, which is the reason for our terminology. 
We note that the family of CNS-graphs constitutes a subclass of closed neighborhood distinct graphs, 
also known as point-determining graphs which were firstly defined and studied by Sumner~\cite{PS}.

\begin{theorem}\label{thm:wc-cns}
If $G$ is a well-covered $(C_4,C_5,C_7)$-free graph, then $G$ is codismantlable. In particular,
such a graph is vertex-decomposable.
\end{theorem}
\begin{proof}
Assume that $G$ is a minimal counterexample to our claim, that is, $G$ is a well-covered $(C_4,C_5,C_7)$-free graph without
any codominated vertex, i.e., it is a CNS-graph, while having the least possible number of vertices with these properties.  
Let $x\in V$ be a non-isolated vertex. Then by Lemma~\ref{lem:unmixed}, the graph
$H:=G-N_G[x]$ is a well-covered graph which is also $(C_4,C_5,C_7)$-free. Therefore, $H$ can not be a CNS-graph by the minimality
of $G$. So, it has a codominated vertex, say $y\in V(H)$. Thus, $H$ is codismantlable to $H-y$. But, $H-y$ is again
a well-covered $(C_4,C_5,C_7)$-free graph; hence, it is either an edgeless graph or else it 
has a codominated vertex. Continuing in this way, we therefore conclude that $H$ is
a codismantlable graph. Then $G$ must be codismantlable by Lemma~\ref{lem:lnk-codis}, a contradiction. 

Finally, since $G$ is codismantlable, there exists a codominated vertex, say $x'\in V(G)$. It follows that both
$G-x'$ and $G-N_G[x']$ are well-covered  $(C_4,C_5,C_7)$-free graphs. By induction both are vertex-decomposable, hence, so is $G$. 
\end{proof}

We note that Theorem~\ref{thm:wc-cns} ensures that the vertex-decomposability, 
codismantlability and Cohen-Macaulayness are all equivalent for a well-covered $(C_4,C_5,C_7)$-free graph.

An easy consequence of Theorem~\ref{thm:wc-cns} is the fact that CNS-graphs lack the Cohen-Macaulay property in general. 
We should recall that only well-covered cycles are of length $3,4,5$ 
and $7$, and among these $C_3$ and $C_5$ are the ones that are also Cohen-Macaulay~\cite{FVT1}. 

\begin{corollary}\label{cor:cns-ncm}
If $H$ is a $(C_4,C_5,C_7)$-free CNS-graph, then $H$ is not well-covered. In particular, such a graph
is not Cohen-Macaulay.
\end{corollary}

Note that if $x$ is a shedding vertex of a bipartite graph, it has to be codominated by Proposition~\ref{prop:bip-vd}, 
and if $y\in V\backslash \{x\}$ is the vertex for which $N_G[y]\subseteq N_G[x]$, then $y$ is of degree $1$, 
since $G$ is triangle-free. This observation provides a necessary condition
on bipartite vertex-decomposable graphs (compare with Theorem 1.1 of~\cite{VTV} and Theorem 2.10 of~\cite{VT}).

\begin{corollary}\label{cor:bip-vd}
If a bipartite graph $G$ is vertex-decomposable, then it is codismantlable.
\end{corollary}
\begin{proof}
By Proposition~\ref{prop:bip-vd}, a shedding vertex, say $x\in V$ of $G$ must be codominated, that is,
$G$ is codismantlable to $G-x$, so the result follows by an induction.
\end{proof}

It should be noted by recalling the graph $G$ in Figure~\ref{contr-1} that
the converse of Corollary~\ref{cor:bip-vd} does not hold in general, unless $G$ is well-covered.

\begin{figure}[ht]
\begin{center}
\psfrag{A}{$C_7$}\psfrag{B}{$P_{10}$}\psfrag{C}{$P_{13}$}\psfrag{D}{$P_{14}$}\psfrag{E}{$Q_{13}$}
\psfrag{x}{$x$}\psfrag{a}{$a$}\psfrag{b}{$b$}\psfrag{s}{$s$}\psfrag{t}{$t$}\psfrag{u}{$u$}
\includegraphics[width=4.8in,height=2.9in]{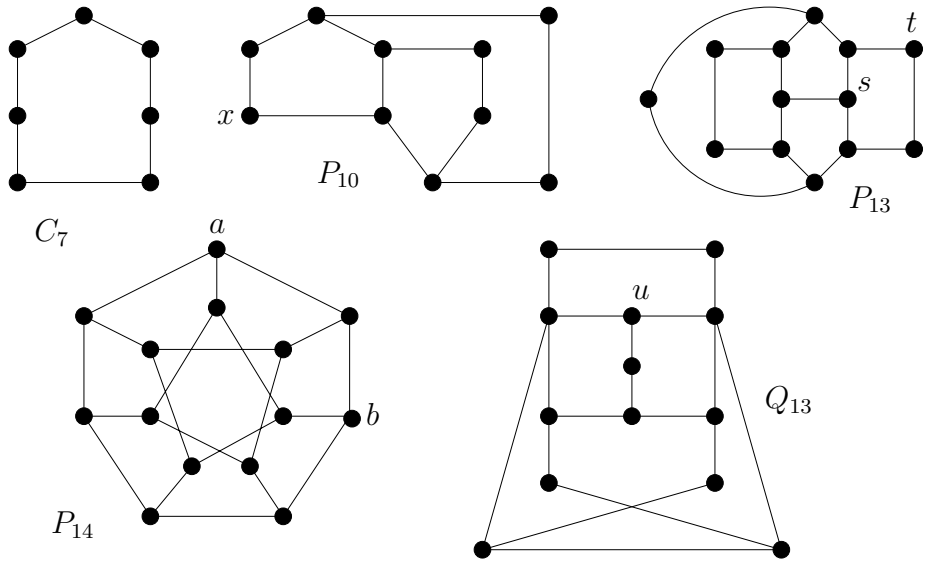}
\end{center}
\caption{Orphans in well-covered graphs with girth at least five.}
\label{orphan}
\end{figure}

Our final aim in this section is to prove the equivalence of vertex-decomposability and Cohen-Macaulayness for
well-covered graphs of girth at least five. We begin with the characterization of well-covered graphs of girth at least five
which is due to Finbow, Hartnell and Nowakowski~\cite{FHN}. They show that any such graph either admits a particular
partition of its vertices or else it is an \emph{orphan}, that is, it is isomorphic to one of five (nontrivial) exceptional 
graphs (see Figure~\ref{orphan}).  
A $5$-cycle in a graph $G$ is said to be \emph{basic} if
it contains no adjacent vertices of degree greater than or equal to $3$ in $G$. Let $\PC$ be the set of graphs $G$ such that
vertex set of $G$ can be partitioned into two disjoint subsets $P$ and $C$ such that $P$ contains the vertices incident with pendant
edges and the pendant edges form a perfect matching of $P$, and $C$ contains the vertices of the basic $5$-cycles and the basic
$5$-cycles form a partition of $C$. 
We recall that when a graph $G$ is well-covered, then the notions of shedding and extendable vertices coincide in $G$. 

\begin{theorem}\cite{FHN}\label{thm:pg}
Let $G$ be a connected well-covered graph with $\gi(G)\geq 5$. The graph $G$ belongs to $\PC$ if and only if
$G$ contains a shedding vertex.
\end{theorem}  

\begin{theorem}\cite{FHN}\label{thm:pgt}
Let $G$ be a connected well-covered graph with $\gi(G)\geq 5$. Then $G$ is in $\PC$ or $G$ is isomorphic
to one of $K_1$, $C_7$, $P_{10}$, $P_{13}$, $Q_{13}$ or $P_{14}$.
\end{theorem}  

We first verify that (nontrivial) orphans can not be Cohen-Macaulay. 
\begin{proposition}\label{prop:orphans}
The graphs $C_7$, $P_{10}$, $P_{13}$, $Q_{13}$ or $P_{14}$ are not Cohen-Macaulay.
\end{proposition}
\begin{proof}
We already know that $C_7$ is not Cohen-Macaulay.
Following the labeling in Figure~\ref{orphan}, it is clear that $P_{10}-N_{P_{10}}[x]\cong C_7$,
$P_{13}-N_{P_{13}}[\{s,t\}]\cong C_7$, $P_{14}-N_{P_{14}}[\{a,b\}]\cong C_7$ and
$Q_{13}-N_{Q_{13}}[u]\cong C_7\cup K_2$. Since the links in a Cohen-Macaulay graphs are Cohen-Macaulay,
and if a graph is Cohen-Macaulay, then so are its connected components, we therefore conclude that
none of these graphs can be Cohen-Macaulay.
\end{proof}

\begin{theorem}\label{thm:gi5}
Let $G$ be a well-covered graph with $\gi(G)\geq 5$ and $E(G)\neq \emptyset$. Then the followings are equivalent:
\begin{itemize}
\item[(1)] $G$ is vertex-decomposable,
\item[(2)] $G$ is Cohen-Macaulay,
\item[(3)] $G\in \PC$.
\end{itemize}
\end{theorem}
\begin{proof}
The implication $(2)\Rightarrow (3)$ is Proposition~\ref{prop:orphans}, while the implication
$(1)\Rightarrow (2)$ is well-known to hold and graph $G$. It remains to prove $(3)\Rightarrow (1)$.

So, let $G\in \PC$. We may assume without loss of generality that $G$ is connected. 
This in particular implies that $G$ has a shedding vertex, say $x\in V$ by Theorem~\ref{thm:pg}. For such a vertex, the
graphs $G-x$ and $G-N_G[x]$ are both well-covered by Lemma~\ref{lem:unmixed}, and both have girth at least five. 
If $G-x$ and $G-N_G[x]$ are both edgeless, then they are vertex-decomposable, hence, so is $G$.
Suppose that each of these graphs contains at least one edge. We claim that $G-x, G-N_G[x]\in \PC$.
In fact, $C_7$ can not be a connected component of $G-N_G[x]$, since otherwise we would have $G\notin \PC$. On the other hand,
the $5$-cycles contained in any of the other orphans are not basic so that none of these graphs can be a connected component of $G-N_G[x]$;
hence, $G-N_G[x]\in \PC$. By a 
similar argument, we have $G-x\in \PC$ (up to isolated vertices). 
We may therefore apply induction on the order of $G$ to deduce that $G-x$ and $G-N_G[x]$ are both vertex-decomposable, so $G$ is as well.
\end{proof}
\begin{remark}
We note that Theorem~\ref{thm:gi5} is also discovered independently by Hoang, Minh and Trung~\cite{HMT} (compare to Theorem 3.1 in~\cite{HMT}). 
\end{remark}

\section{Regularity of Codismantlable graphs}

In this section, we study regularity of codismantlable graphs.
Consider the minimal free graded resolution of $R/I(G)$ as an $R$-module:
\begin{equation*}
0\to \bigoplus_j R(-j)^{\beta_{i,j}}\to \ldots \to \bigoplus_j R(-j)^{\beta_{1,j}}\to R\to R/I(G)\to 0.
\end{equation*}
The \emph{Castelnuovo-Mumford regularity} or simply the \emph{regularity} of $R/I(G)$ is defined as
\begin{equation*}
\reg(G):=\reg(R/I(G))=\max \{j-i\colon \beta_{i,j}\neq 0\}.
\end{equation*}

Most of our results will be based on the following inductive machinery on computing the regularity of certain graphs.

\begin{theorem}\cite{MV}\label{thm:mv}
Let $\F$ be a family of graphs containing any edgeless graph and let $\eta\colon \F\to \N$
be a function satisfying $\eta(G)=0$ for any edgeless graph $G$, and such that given $G\in \F$
with $E(G)\neq \emptyset$, there exists $x\in V(G)$ such that the following two conditions hold:
\begin{itemize}
\item[(i)] $G-x$ and $G-N_G[x]$ are in $\F$,\\
\item[(ii)] $\eta(G-N_G[x])<\eta(G)$ and $\eta(G-x)\leq \eta(G)$.
\end{itemize}
Then $\reg(G)\leq \eta(G)$ for any $G\in \F$.
\end{theorem}

\begin{lemma}\label{lem:reg-codom}
If $x$ is a codominated vertex of a graph $G$, then $\im(G-x)\leq \im(G)$ and $\im(G-N_G[x])<\im (G)$.
\end{lemma}
\begin{proof}
The first inequality trivially holds, so we verify the second. Suppose that $M\subseteq E(G-N_G[x])$
is an induced matching. Since $x$ is codominated in $G$,
there exists a vertex $y\in V$ satisfying $N_G[y]\subseteq N_G[x]$. It then follows that $M\cup \{xy\}$
is an induced matching in $G$ so that $\im(G-N_G[x])<\im (G)$ as required. 
\end{proof}

\begin{theorem}\label{thm:unmixed-reg}
If $G$ is a $(C_4,C_5)$-free vertex-decomposable graph, then $\reg(G)=\im(G)$.
\end{theorem}
\begin{proof}
By a result of Katzman~\cite{MK}, we already know that $\im(G)\leq \reg(G)$ for any graph $G$.
So, it is suffices to show that $\im(G)$ is an upper bound when $G$ fulfills the required condition.

Suppose that $G$ is a $(C_4,C_5)$-free vertex-decomposable graph. Since $G$ is vertex-decomposable,
it must have a shedding vertex $v\in V$, and such a vertex is codominated by Theorem~\ref{thm:shedding}.
Note also that $G-v$ and $G-N_G[v]$ are $(C_4,C_5)$-free vertex-decomposable graphs; hence, 
Theorem~\ref{thm:mv} applies. The result follows by Lemma~\ref{lem:reg-codom}.
\end{proof}

We first state an easy consequence of Theorem~\ref{thm:unmixed-reg} which is originally due to H\'a and Van Tuyl~\cite{HVT} 
by recalling that chordal graphs are vertex-decomposable~\cite{RW1}.

\begin{corollary}
If $G$ is a chordal graph, then $\reg(G)=\im(G)$.
\end{corollary}

Let $G$ be a graph, and $\HE$ be a family of graphs. The $\HE$-{\it cover number} of $G$ is
the minimum number of subgraphs $H_1,\ldots,H_r$ of $G$ such that every $H_i\in \HE$ and $\cup E(H_i)=E(G)$.
In particular, we denote by $\cd(G)$, the {\it cochordal cover number} of $G$. 

Woodroofe~\cite{RW2} has recently shown that the cochordal cover number of a graph provides an upper bound for its regularity.
He also provides examples of graphs whose cochordal covering numbers
could be far from their regularities. 

The graph $K_{1,m}$ is called a {\it star graph}, and the graph obtained by joining the center
vertices of two stars by an edge is known as a {\it double-star}.
A subset $F\subseteq E$ is called an {\it edge-dominating set} of $G$ if each edge of $G$ either belongs to $F$ or is
incident to some edge in $F$, and the {\it edge-domination number} $\gamma'(G)$ of $G$ is the minimum cardinality
of an edge-dominating set of $G$. We remark that $\gamma'(G)$ equals to the minimum size of a maximal matching
for any graph $G$~\cite{YG}.
 
\begin{lemma}\label{lem:double-star}
If $H$ is a connected cochordal graph with $\gi(H)\geq 5$, then $H$ is a star or a double-star.
\end{lemma}
\begin{proof}
Let $H$ be a cochordal graph. Note that $H$ is $P_5$-free, since otherwise the complement of $H$ contains an
induced $C_4$, which is not possible. Therefore, the diameter of $H$ is at most three. If $\di(H)=2$, then $H$
is a star. On the other hand, if $\di(H)=3$, then $H$ has an induced $P_4=(\{v_1,v_2,v_3,v_4\},\{v_1v_2,v_2v_3,v_3v_4\})$.
We claim that $H$ is a double-star with the center-edge $v_2v_3$. In fact, 
if $a\in N_H(v_2)\backslash \{v_3\}$ and $b\in N_H(v_3)\backslash \{v_2\}$, then we must have $ab\notin E(H)$,
since $\gi(H)\geq 5$. By a similar reasoning, we have $uw\notin E(H)$ whenever $u,w\in N_H(v_i)$ for $i=2,3$ so that
the claim follows.
\end{proof}
 
When $H$ is a cochordal subgraph of a graph $G$ with $\gi(G)\geq 5$, we now know from Lemma~\ref{lem:double-star}
that $H$ is either a star or a double-star. If $H$ is a star with a center vertex, say $v$, 
we extend it to a double star by choosing a neighbor $u\in N_G(v)$, and adding every edge (if any) to $H$ 
whose one end is $u$. We then call the edge $e=uv$ as the \emph{center-edge} of $H$. We may therefore consider a star
in $G$ to be a (degenerate) double-star.

\begin{theorem}\label{thm:cd-dom}
If $G$ is a connected graph with $\gi(G)\geq 5$, then $\cd(G)=\gamma'(G)$.
\end{theorem}
\begin{proof}
If $G$ is isomorphic to $K_2$, the statement is trivial. Suppose that $\HE=\{H_1,\ldots,H_k\}$ is a cochordal edge cover of $G$, and let
$f$ be an edge of $G$. Since $\HE$ is an edge cover, there exists an $i\in [k]$ such that $f\in E(H_i)$. However, $H_i$ is a star
or a double-star by Lemma~\ref{lem:double-star} so that the edge $f$ is either the center edge of $H_i$ or it is incident to the center
edge of $H_i$. Therefore, the set of center edges of $H_i$'s forms an edge-dominating set of $G$.

For the converse, let $\{e_1,\ldots,e_k\}$ be an edge-dominating set. If we define $R_i$ to be the subgraph consisting of
those edges incident to the edge $e_i$ for any $1\leq i\leq k$, then each $R_i$ is a double-star, and the family
$\{R_1,\ldots,R_k\}$ is clearly a cochordal cover for $G$. 
\end{proof}

Our final aim in this section is to understand the structure of graphs satisfying $\im(G)=\m(G)$. 
When $\im(G)=\m(G)>1$ for a graph $G$, Kobler and Rotics~\cite{KR} have a nice description of how the largest induced 
matchings arise in $G$.
 
\begin{lemma}\cite{KR}\label{lem:KR}
Let $G$ be a connected graph such that $\im(G)=\m(G)>1$. Then any edge
in a maximum induced matching in $G$ is either a pendant or a triangle edge.
\end{lemma}

We should note that Lemma~\ref{lem:KR} remains valid even if $G$ is disconnected.

\begin{proposition}\label{prop:im-m}
Let $G=(V,E)$ be a graph such that $\im(G)=\m(G)$. If $x$ is a vertex of $G$ for which $N_G(y)=\{x\}$
for some vertex $y\in V$, then $\im(G-x)=\m(G-x)=\im(G)-1$.
\end{proposition}
\begin{proof}
It suffices to show that 
\begin{equation*}
\m(G)-1\geq \m(G-x)\geq \im(G-x)\geq \im(G)-1.
\end{equation*}
The rightmost inequality follows because an induced matching of $G-x$ is an induced matching of $G$, and at most
one edge in an induced matching of $G$ can contain $x$. The middle inequality is because an induced matching is a matching.
Finally, the edge $xy$ can be added to any matching in $G-x$, giving the leftmost inequality. 
\end{proof}

\begin{theorem}\label{thm:tf-r-i-m}
Let $G$ be a graph. Then $G$ is codismantlable if $\im(G)=\m(G)$.
\end{theorem}
\begin{proof}
If $\im(G)=\m(G)=0$, then $G$ is an edgeless graph so that there is nothing to prove. 
On the other hand, if $\im(G)=\m(G)=1$, then either $G$ is a star or else it is isomorphic to $K_3$ (up to isolated vertices);
hence, $G$ is codismantlable in either case.

We let $\im(G)=\m(G)>1$. By Lemma~\ref{lem:KR}, any edge in a largest induced matching is either a pendant or a triangle edge of $G$.
If $e=uv\in E$ is a triangle edge, then $u$ and $v$ are of degree $2$ and have a common neighbor, say $w\in V$.
Now, the vertex $w$ is codominated in $G$, since $N_G[u]=N_G[v]\subseteq N_G[w]$, that is,
$G$ is codismantlable to $G-w$. Note also that $\im(G-w)=\im(G)=\m(G)=\m(G-w)$ for such a vertex.
If $f=xy\in E$ is a pendant edge with $N_G(y)=\{x\}$, the vertex $x$ is codominated in $G$,
so $G$ is codismantlable to $G-x$. By Proposition~\ref{prop:im-m}, we have $\im(G-x)=\m(G-x)$. The theorem
follows by induction.
\end{proof}

\begin{figure}[ht]
\begin{center}
\includegraphics[width=1.6in,height=1.5in]{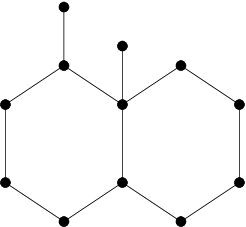}
\end{center}
\caption{A codismantlable graph $H$ with $\reg(H)=\im(H)<\m(H)$.}
\label{Fig:im-m}
\end{figure}

The converse of Theorem~\ref{thm:tf-r-i-m} does not hold in general. For the codismantlable graph $H$ depicted in 
Figure~\ref{Fig:im-m}, we have $\reg(H)=\im(H)=3$ and $\m(H)=6$.

\section{Orientations and vertex-decomposable graphs}
In this section, we introduce an operation on digraphs that allows us to construct a vertex-decomposable graph whenever
the source digraph is acyclic. We recall that when $\ve{G}=(V,E)$ is a digraph, that is, its set of edges $E$ consists of 
ordered pairs of vertices, the digraph $\ve{G}$ is said to be {\it acyclic}, if there exists no directed cycles in $\ve{G}$.

Let $\ve{G}=(V,E)$ be a digraph. We write $x\Rightarrow y$ when there exists a dipath (directed path) in $\ve{G}$ 
starting from $x$ and ending at $y$, and define the \emph{enemy-set} of $u\in V$ by $A(u):=\{v\in V\colon v\Rightarrow u\}$ and
set $A[u]:=A(u)\cup \{u\}$ for each $u\in V$.

\begin{definition}
For a given digraph $\ve{G}=(V,E)$, we define its \emph{common-enemy graph} $CE(\ve{G})$ to be the (simple) graph on $V$ such that
$xy\in E(CE(\ve{G}))$ if and only if $x\neq y$ and $A[x]\cap A[y]\neq \emptyset$. In particular, we call the independence complex 
$I(CE(\ve{G}))$ of $CE(\ve{G})$ as the \emph{common-enemy complex} of $\ve{G}$.
\end{definition}

\begin{remark}
The notion of the common-enemy graph of an acyclic digraph is firstly defined by Cohen (see~\cite{RL} for a survey article),
slightly differently how we use it here. The one we consider corresponds to the usual common-enemy graph taken
on the transitive completion of $\ve{G}$.
\end{remark}

\begin{theorem}\label{thm:cf-acyclic}
If $\ve{G}$ is an acyclic digraph, then the graph $CE(\ve{G})$ is vertex-decomposable and codismantlable.
\end{theorem}
\begin{proof}
In order to simplify the notation, we will write $H:=CE(\ve{G})$.
Since $\ve{G}$ is acyclic, it must contain a vertex $x$ of out-degree zero,
and without loss of generality, we may assume that it has a positive in-degree. We claim that such a vertex
is codominated in $H$. Let $y\in V\backslash \{x\}$ be a vertex of $\ve{G}$ such that $(y,x)\in E(\ve{G})$.
Then we need to have $A[y]\subseteq A[x]$; hence, $N_H[y]\subseteq N_H[x]$. Moreover,
we have $H-x\cong CE(\ve{G}-x)$ and $H-N_H[x]\cong CE(\ve{G}-N_H[x])$. In fact, $CE(\ve{G}-x)$ and 
$CE(\ve{G}-N_H[x])$ are clearly subgraphs of $H-x$ and $H-N_H[x]$ respectively.
If $ab\in E(H-x)$, then there exists $z\in V$ satisfying $a\Leftarrow z \Rightarrow b$ 
for some $z\in V$. Note that $z\neq x$, since $x$ has out-degree zero. Therefore, $ab\in E(CE(\ve{G}-x))$.

Similarly,  assume that $cd\in E(H-N_H[x])$ so that $c\Leftarrow u \Rightarrow d$ for some vertex
$u\in V\backslash \{x\}$. Let $u=x_0,x_1,\ldots, x_k=c$ be the vertices of the dipath from $u$ to $c$.
If $x_i\in N_H[x]$ for some $0\leq i\leq k$, then there exists $w\in V$ such that 
$x_i\Leftarrow w \Rightarrow x$. However, this implies $c\in N_H[x]$, since $w\in A[c]\cap A[x]$, which is not possible.

Therefore, the vertex $x$ is a shedding vertex of $H$ for which $H-x$ and $H-N_H[x]$ are common-enemy graphs of
acyclic digraphs; hence, the vertex-decomposability and codismantlability of $CE(\ve{G})$ both follow by an induction.
\end{proof}

Acyclic digraphs naturally arise when we consider cover digraphs of partially ordered sets (poset for short).
We recall that if $P=(X,\leq)$ is a (finite) poset with $x,y\in X$ such that $x\neq y$, then $x$ covers $y$ in $P$,
if $y\leq x$ and there exists no $z\in X\backslash \{x,y\}$ satisfying $y<z<x$.
Then, the cover digraph $\cov(P)$  of $P$ is defined to be the digraph on $X$ such that
$xy\in E(\cov(P))$ if and only if $x\neq y$ and $x$ covers $y$ in $P$.
In the context of cover digraphs, the notion of common-enemy graphs corresponds to a well-known class of graphs. 
In detail, the upper-bound graph $UB(P)$ of $P$ is defined to be the graph on $X$ with
$xy\in E(UB(P))$, where $x\neq y$ and there exists $z\in X$ such that $x,y\leq z$ (see~\cite{MZ}).
Therefore, the following is a consequence of Theorem~\ref{thm:cf-acyclic}.

\begin{corollary}\label{cor:ub-vd}
The upper-bound graph of any poset is vertex-decomposable and codismantlable.
\end{corollary}

Our final aim is to provide an answer to a question raised in~\cite{MMCRTY} regarding the existence of sequentially 
Cohen-Macaulay graphs
carrying some specific properties. Having this in mind, we first introduce an operator ``the edge-clique-whiskering''
that turns any graph into a vertex-decomposable graph. Beforehand, we should note that Cook and Nagel~\cite{CN}
have also introduced such an operation for a different purpose. They consider
vertex-disjoint clique partitions in their construction, while we work with edge-disjoint such partitions,
so the one considered here applies in somewhat different circumstances from theirs.

\begin{definition}
Let $G=(V,E)$ be a connected graph with $E\neq \emptyset$, and let $\Pi:=\{W_1,\ldots,W_k\}$ be a family of cliques of $G$, each having
order at least $2$. We call $\Pi$ as an \emph{edge-clique partition} of $G$, if $\{E(G[W_1]),\ldots, E(G[W_k])\}$ is a partition of $E$.
We then define the \emph{edge-clique-whisker} $G^{\Pi}$ of $G$ with respect to the partition $\Pi$ to be the graph
on $V(G^{\Pi}):=V\cup \{x_1,\ldots,x_k\}$ with $E(G^{\Pi}):=E\cup \{ux_i\colon u\in V(W_i), i\in [k]\}$.
We particularly call the newly added vertex $x_i$ as the \emph{cone-vertex} of $G^{\Pi}$ corresponding
to the clique $W_i$.
\end{definition}

We note that every connected graph admits an edge-clique partition, for instance, its set of edges forms such a partition.

\begin{theorem}\label{thm:c-w-vd}
The graph $G^{\Pi}$ is vertex-decomposable and codismantlable for any edge-clique partition $\Pi$ of $G$.
\end{theorem}
\begin{proof}
We consider the digraph $\ve{G}_{\Pi}=V(G^{\Pi})$ with directed edges $(x_i,u)\in E(\ve{G}_{\Pi})$ if
$u\in W_i$ for any $i\in [k]$ and no more edges. It is clear that the common-enemy graph $CE(\ve{G}_{\Pi})$
is isomorphic to $G^{\Pi}$. Furthermore, $\ve{G}_{\Pi}$ is acyclic, since there are edges only among the vertices of the source
graph $G$ and corresponding cone-vertices. Therefore, our claim follows from Theorem~\ref{thm:cf-acyclic}.
\end{proof}

\begin{figure}[ht]
\begin{center}
\includegraphics[width=3.1in,height=1.5in]{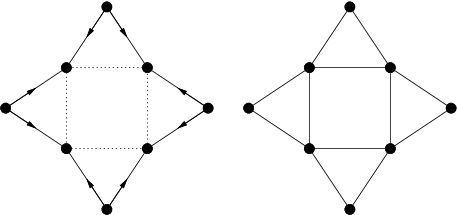}
\end{center}
\caption{The edge-clique-whiskering of the $4$-cycle.}
\label{or-2}
\end{figure}

When Theorems~\ref{thm:unmixed-reg} and \ref{thm:c-w-vd} are combined, we obtain the following. 
\begin{corollary}\label{cor:reg-whisker}
If $G$ is a $(C_4,C_5)$-free graph, then $\reg(G^{\Pi})=\im(G^{\Pi})$ for any edge-clique partition $\Pi$ of $G$.
\end{corollary}

As we mentioned earlier, the following question is raised in~\cite{MMCRTY}. 
\begin{question}\cite{MMCRTY}
Let $G$ be a SCM graph with $2n$ vertices which are not isolated such that the vertex cover number of $G$ equals to $n$. 
Then do we have the following statements?
\begin{itemize}
\item[$(1)$] $G$ is vertex-decomposable,
\item[$(2)$] $\reg(G)=\im(G)$.
\end{itemize}
\end{question}

In fact, it was originally asked whether such a graph admits a vertex of degree-one while \cite{MMCRTY} was a preprint.
Now, following Theorem~\ref{thm:c-w-vd}, the edge-clique-whiskering $C_n^{\Pi}$ of the cycle $C_n$ is vertex-decomposable; 
hence it is sequentially Cohen-Macaulay with the required properties without any vertex of degree-one 
for any $n\geq 4$, where $\Pi=E(C_n)$ (see Figure~\ref{or-2} for $n=4$). 
Furthermore, when $n\geq 6$, we have $\reg(C_n^{\Pi})=\im(C_n^{\Pi})$
by Corollary~\ref{cor:reg-whisker}.

\section{conclusion}
Our results combine algebraic and topological notions such as (sequentially) Cohen-Macaulayness and 
Castelnuovo-Mumford regularity via graph theory by codismantlable graphs. The family of codismantlable graphs
seems to be an interesting graph class in their own right, and many intriguing questions remain to be answered 
concerning the role of codismantling operation in combinatorics. We hope to return to these in future.
Meantime, we will discuss some additional problems related to our present work.

\begin{problem}\label{prob:wccodis-vd}
Is there any well-covered codismantlable graph which is not vertex decomposable?
\end{problem}

At this point, it is also reasonable to ask under what conditions on a codismantlable graph $G$, we may guarantee that 
$\reg(G)=\im(G)$.

\begin{problem}\label{prob:cns-cm}
Is there any $(C_4,C_5)$-free CNS-graph which is Cohen-Macaulay over some field?
\end{problem}
We believe that a negative answer to Problem~\ref{prob:cns-cm} would be of no surprise. In detail,
we call a connected induced subgraph $H$ of a graph $G$ \emph{separated} in $G$, if there exists an independent set
$U$ of $G$ such that $H$ is a connected component of $G-N_G[U]$. Note that in a  Cohen-Macaulay graph,
no cycle of length at least $6$ can be separated. So, one way to attack Problem~\ref{prob:cns-cm} is
to check whether such a cycle can be separated in a well-covered, $(C_4,C_5)$-free CNS-graph.

\section*{Acknowledgments}
We would like to thank the anonymous referee for the very useful comments and
detailed corrections, which we found very constructive and helpful to improve our manuscript.
We have benefited greatly from his/her expertise in the field. In particular, we are grateful to him/her for suggesting 
the current short proof of Proposition~\ref{prop:im-m} which has replaced the original one.


\end{document}